\RequirePackage{fix-cm}
\documentclass[smallextended]{svjour3}       
\smartqed  
\usepackage{graphicx}
\usepackage{latexsym}

\usepackage{amsfonts,amsmath}

\usepackage{amssymb}

\usepackage[mathscr]{eucal}

\begin{document}

\title{
Schr\"odinger dynamics and optimal transport of measures on the  torus
}


\author{Lorenzo Zanelli}


\institute{L. Zanelli \at
              University of Padova \\
              Department of Mathematics Tullio Levi-Civita\\
              \email{lorenzo.zanelli@unipd.it}  }

\date{Received: date / Accepted: date}

\maketitle

\begin{abstract}
The aim of this paper is to recover displacement interpolations of probability measures, in the sense of the Optimal Transport theory, by semiclassical measures associated with solutions of  Schr\"odinger's equations defined on the flat torus. Under an additional assumption, we show the completing viewpoint by proving that a family of displacement interpolations can always be viewed as these time dependent semiclassical measures.    
\end{abstract}

\section{Introduction}
\label{intro}
Let $\Bbb T^n := ( \Bbb R / 2 \pi \Bbb Z )^n$, $V \in C^\infty (\Bbb T^n)$.  Let us consider the  Schr\"odinger equation 
\begin{equation}
\label{Sch}
i \hbar \partial_t \psi_\hbar (t,x) =  - \frac{\hbar^2}{2m} \Delta_x \psi_\hbar (t,x) + V(x) \psi_\hbar (t,x). 
\end{equation}
The Schr\"odinger dynamics can be given by the one parameter group of unitary operators $U_\hbar (t) := e^{- i \widehat{H} t / \hbar}$, $\widehat{H} := - \hbar^2 \Delta_x  / 2m + V(x)$,  acting on $L^2 (\Bbb T^n)$.
Thus, we consider the solution of (\ref{Sch}) as $\psi (t,x) := (U_\hbar (t)  \varphi_\hbar) (x)$ for 
initial data in a class of WKB - wave functions 
\begin{equation}
\label{def-wave}
\varphi_\hbar (x) = a_\hbar (x) \, e^{i S_+  (x) / \hbar},
\end{equation}
the related semiclassical probability measures $\omega_t \in \mathcal{P} (\Bbb T^n \times \Bbb R^n)$ associated with the path $\psi_\hbar (t,\cdot \,)|_{0 \le t \le 1}$ and study $\sigma_t := \pi_\sharp \, \omega_t  \in \mathcal{P} (\Bbb T^n)$. In order to select (\ref{def-wave}) we fix the phase as a Lipschitz continuous weak KAM solution of positive type for the following Hamilton-Jacobi equation (see \cite{F})
\begin{equation}
\label{H-J00}
\frac{1}{2m} |\nabla_x S_+ (x)|^2 + V(x) = \max_{y \in \Bbb T^n} V(y).
\end{equation}
The reason why we select such solutions is the time forward invariance property of the graph under the Hamiltonian flow of  $H := |p|^2 / 2m + V(x)$, namely $\phi_{H}^t ( {\rm Graph}( \nabla_x S_+ ) ) \subseteq {\rm Graph}( \nabla_x S_+ )$ $\forall t \ge 0$, as shown by  Thm 4.9.3 in \cite{F}.\\  
The amplitudes in (\ref{def-wave}) are selected by  $a_\hbar \in H^1 (\Bbb T^n;\Bbb R)$ where $\| a_\hbar \|_{L^2} = 1$, $\| \hbar \nabla a_\hbar \|_{L^2} \to 0$ as $\hbar \to 0$, $a^2 (x) dx \rightharpoonup \sigma_0$ weakly as  measures on $\Bbb T^n$, and such that ${\rm supp}(\sigma_0) \subseteq {\rm dom}(\nabla S_+)$, $\sigma_0 \in \mathcal{P}_{ac} (\Bbb T^n)$, i.e. Borel probability measures which are absolutely continuous with respect to Lebesgue. As we will see, the assumption of the absolute continuity of $\sigma_0$ with respect to the Lebesgue measure $\mathcal{L}^n$ turns out to be useful in order to make a full relationship between a class of optimal transport problems of measures on $\Bbb T^n$ and the semiclassical measures arising from our Schr\"odinger equation (\ref{Sch}).
 
The first result of our paper show that the above family of projected semiclassical measures $(\sigma_t )_{0 \le t \le 1}$ is a displacement interpolation between the Borel probability measures  $\sigma_0, \sigma_1 \in \mathcal{P} (\Bbb T^n)$ in the sense of  Optimal Transport theory (see Thm 7.21 in \cite{V}). More precisely, here we deal with the minimum curves for the Action functional 
\begin{equation} 
\label{cost00}
\inf_\gamma \Big(  \int_{\Omega}  \int_0^1 \frac{m}{2} |\dot{\gamma} (t,\zeta) |^2 - V(\mathbf{\gamma} (t,\zeta))  dt \, {\rm d}\mathbb{P} (\zeta) \Big)
\end{equation}
where the infimum is over all the random curves $\gamma : [0,1] \times \Omega \longrightarrow \Bbb T^n$ such that ${\rm Law} (\gamma (t, \cdot \,) ) = \sigma_t$.  In particular, we are interested to deal with the family of the above displacement interpolations coming from  solutions $\sigma \in C([0,1];  \mathcal{P} (\Bbb T^n) )$ of the continuity equation  in the measure sense
\begin{equation}
\label{eqTR01}
\partial_t \sigma_t (x) +  {\rm div}_x  \Big( \frac{1}{m} \nabla_x S_+  (x)  \sigma_t (x) \Big) = 0
\end{equation}
for arbitrary fixed $\sigma_0 \in  \mathcal{P}_{ac} (\Bbb T^n)$. \\  
Before to state precisely the main results of the paper, we underline that they are mainly based on some meaningful arguments of semiclassical Analysis and Optimal Transport theory. The first one is that 
the continuous paths of semiclassical measures $\omega_t$ associated to the solution of the Schr\"odinger equation solve the Liouville equation in the measure sense,  
\begin{equation}
\label{eqTR02}
\partial_t \, \omega_t (x,p) +  p \cdot \nabla_x  \, \omega_t (x,p) - \nabla_x V (x) \cdot  \nabla_p  \, \omega_t (x,p) = 0
\end{equation}
as firstly shown in  \cite{L-P}  within the euclidean setting (and many others under various assumptions, see \cite{A-F-P} and the references therein) and recently in \cite{TP-Z} within the toroidal setting.  
The second main ingredient is that all the semiclassical measures of $\varphi_\hbar$ as in (\ref{def-wave}) take the form
\begin{equation} 
\label{omega00}
\omega_0 (x,p) = \delta ( p - \nabla_x S_+ (x) ) \sigma_0 (x).
\end{equation}
This property becomes meaningful in view of  the time forward invariance property of ${\rm Graph}( \nabla_x S_+ )$ under the Hamiltonian flow, as proved in \cite{F}.  Furthermore, we take into account the equations linked to displacement interpolations  of measures, as described in \cite{V}, exhibiting in our paper the simple form (\ref{eqTR01}). Finally, we take into account the results  on the existence of the transport maps ${\rm T}_t : \mathcal{P}_{ac} (\Bbb T^n) \to \mathcal{P} (\Bbb T^n)$  which  solve the Monge problem (see \cite{F-F}, \cite{K1}, \cite{K2} and the references therein) for the cost function $c^{0,t} (x,y) := \inf_\gamma \int_0^t L(\gamma,\dot{\gamma}) d\tau$, $L := m |\xi|^2 /2 - V(x)$, $\gamma$ are at least $C^1$ and fulfill $\gamma (0) = x$, $\gamma (t) = y$,  and they provide displacement interpolations by  
\begin{equation} 
\sigma_t = ({\rm T}_t)_\sharp \sigma_0.
\end{equation}				
In fact, it turns out that any of such transport map read ${\rm T}_t = \pi  \circ \phi_{H}^t (x,\nabla_x f (x))$ for some Lipschitz functions $f : \Bbb T^n \to \Bbb R$ which is convex with respect to the cost function $c^{0,1} (x,y)$ namely for some $\bar{f} : \Bbb T^n \to \Bbb R$  it holds $f(x) = \sup_{y \in \Bbb T^n} (\bar{f}(y) - c^{0,1}(x,y))$. 
In our paper, $f = S_+$  as we will easily see in the Remark \ref{REM-convex}. In the paper \cite{B-B}, the convex condition for $S_+$ is shown from a more general viewpoint involving Monge-Kantorovich duality.

We are now ready to provide the first result of the paper

\begin{theorem}
\label{TH1}
Let $\varphi_\hbar$ be as in (\ref{def-wave}) and $\omega_0 \in \mathcal{P} (\Bbb T^n \times \Bbb R^n)$ be an associated semiclassical measure. 
Let $\phi_{H}^t : \Bbb T^n \times \Bbb R^n \to \Bbb T^n \times \Bbb R^n$ be the Hamiltonian flow of $H =  |p|^2 / 2m  + V(x)$. 
Then, the  $\omega_t := (\phi_{H}^t)_\sharp \omega_0$ is a semiclassical measure associated with  $\psi_\hbar (t,\cdot \,)$ and takes the form $\forall$ $0 \le t \le 1$
\begin{equation} 
\label{omega}
\omega_t (x,p) = \delta ( p - \nabla_x S_+ (x) ) \sigma_t (x)
\end{equation}
and the path $( \sigma_t )_{0 \le t \le 1} \in \mathcal{P} (\Bbb T^n)$ equals  for $\mathcal{L}^1$ - a.e. $0 \le t \le 1$  a continuous displacement interpolation  between $\sigma_0$ and $\sigma_1$ in the sense of (\ref{cost00}) and (\ref{eqTR01}). By defining $\Psi^t (x) := \pi \circ \phi_{H}^t (x,\nabla_x S_+ (x))$
it holds  $\sigma_t = (\Psi^t)_\sharp (\sigma_0)$ $\forall 0 \le t \le 1$, i.e.
\begin{equation} 
\label{omega2}
\int_{\Bbb T^n}  g(x) d\sigma_t (x) = \int_{\Bbb T^n} g (\Psi^t (x)) d\sigma_0 (x) \quad \forall g \in C^\infty (\Bbb T^n).
\end{equation}
\end{theorem}    

In the paper \cite{TP-Z}, time propagated semiclassical measures taking the form  $\omega_t (x,p) = \delta ( p -  P - \nabla_x S_\pm (P,x) ) \sigma_t (P,x)$ are studied when  $P \in \ell \Bbb Z^n$ with $\ell >0$,  $\hbar^{-1} \in \ell^{-1} \Bbb N$, $S_\pm$ are weak KAM solutions of positive or negative type for the Hamilton-Jacobi equation
\begin{equation}
\label{HJP}
\frac{1}{2m} |P + \nabla_x S_\pm (P,x)|^2 = \bar{H}(P)
\end{equation}
where $\bar{H}(P) = \sup_x \inf_{v \in C^\infty} \frac{1}{2m} |P + \nabla_x v(x)|^2 + V(x)$ is the so-called effective Hamiltonian (see for example \cite{E1}). Notice that $\bar{H}(0) = \max_{y \in \Bbb T^n} V(y)$ and that for $P = 0$ the equation (\ref{HJP}) becomes (\ref{H-J00}).  
 In particular, any such $\sigma_t$ is absolutely continuous with respect to the projected (on $\Bbb T^n$) $\pi_\sharp (\mu_P)$ where $\mu_P$ are flow invariant and Action-minimizing measure for 
\begin{equation} 
A[\mu] = \int_{\Bbb T^n \times \Bbb R^n} \frac{m}{2} |\xi|^2 - V(x) - P \cdot \xi \ d\mu (x,\xi).
\end{equation}
This setting ensures the possibility to deal with continuity equation (\ref{eqTR01})  for positive or negative times, and in particular to study the propagated densities $g_\pm \in L^1 (\Bbb T^n)$ satisfying  $\sigma_t (P,x) = g_\pm (t,P,x) \pi_\sharp (\mu_P)$.\\ 
We underline that the time propagation under the Hamiltonian flow of measures with a graph structure as in (\ref{omega00}) with different low regularity momentum profiles, have been also recently studied in \cite{bgmp} as an application to the semiclassical limit of quantum propagation of WKB type wave functions.  The semiclassical localization of Schr\"odinger eigenfunctions on the graph of weak KAM solutions has been studied in \cite{Z2}. The related semiclassical measures are invariant under the Hamiltonian flow and whence are not propagating.\\ 
We also recall that in \cite{Z1} the study of time propagated semiclassical measures of type $\omega_t (x,p) = \delta ( p -  P - \nabla_x S_{+} (P,x) ) \sigma_t (P,x)$ in the framework of optimal transport is discussed. By taking $P=0$ it is recovered the form of the measures (\ref{omega}) but  it is still a larger class, since in the current paper we are assuming that $\sigma_0$ is absolutely continuous with respect to the Lebesgue measure. 
Moreover, in \cite{Z1} such measures are obtained by the time propagation of semiclassical limit of the Wigner transform for different initial data wave functions with respect to the ones discussed in the current paper. The more important observation is that in  \cite{Z1} is it proved the analogous result of the Theorem \ref{TH1}  but it is not proved the converse viewpoint, that we are going to prove here in Theorem \ref{TH2} for the smaller class of semiclassical measures (\ref{omega}). Hence, our paper now provides the complete analysis of the  bridge between  the semiclas\-sical measures  supported on weak KAM graphs solving Hamilton-Jacobi equation and ones solving the optimal transport problem on the torus.\\

In the next, we provide the second result of the paper by a complementary viewpoint with respect to Theorem \ref{TH1}.
\begin{theorem}
\label{TH2}
Let $\sigma_0 \in \mathcal{P}_{ac} (\Bbb T^n)$ and assume the uniqueness for solutions $\sigma \in C([0,1];  \mathcal{P} (\Bbb T^n) )$ of (\ref{eqTR01}). Define the lift $\omega_t := \delta ( p - \nabla_x S_+) \sigma_t \in  \mathcal{P} (\Bbb T^n \times \Bbb R^n)$. Then, there exists $\varphi_\hbar$  in the form (\ref{def-wave}) such that $\omega_0$ is the unique linked semiclassical measure. Moreover,  any $\omega_t$ is a semiclassical measure associated with $\psi_\hbar (t,x) := (U_\hbar (t)  \varphi_\hbar) (x)$.  
\end{theorem}    
Notice that here we have assumed the uniqueness for the solution $\sigma$ of the continuity equation  (\ref{eqTR01}) in  $C([0,1];  \mathcal{P} (\Bbb T^n) )$ where $\mathcal{P} (\Bbb T^n)$ is equipped with the L\'evy-Prokhorov distance of probability measures which metrizes the weak convergence. In view of Theorem 3.1 shown in \cite{A2}, such an assumption is equivalent  to the pointwise uniqueness for the solutions of $\dot{\gamma} = \frac{1}{m} \nabla_x S_+ (\gamma)$. In the Lemma \ref{AB}, we provide a solution $\gamma = \pi \circ \phi_{H}^t (x, \nabla_x S_+ (x))$. However, $S_+ : \Bbb T^n \to \Bbb R$ is Lipschitz continuous and $x \to \nabla_x S_+ (x)$ is continuous on its domain, and this low regularity does not guarantees this property. On the other hand, one can assume the additional regularity $\nabla_x S_+  \in W^{1,\infty}_{{\rm loc}} (\Bbb T^n ; \Bbb R^n)$ and apply the Remark 2.1 in \cite{A2} to ensure such a uniqueness, and thus recover the setting for Theorem \ref{TH2}.\\

Before to conclude, we underline a remarkable open problem about  the link between Optimal transport theory and semiclassical Analysis.  More precisely, to prove the existence and related properties for a bigger set of initial data wave functions  $\varphi_\hbar$ taking a more general form than our (\ref{def-wave}) and recovering, in the semiclassical limit, an arbitrary continuous displacement interpolation  as prescribed in (\ref{cost00}) and without the assumption of absolute continuity of the initial measure with respect to Lebesgue.\\

The content of the paper is the following:  in Section \ref{SEC-2} we introduce some preliminaries on Toroidal Pseudodifferential Operators, Weyl quantization and the well posed setting for semiclassical measures on $\Bbb T^n \times \Bbb R^n$. Within the Section \ref{wKAM}  we provide a resume on some central results of the weak KAM theory for Hamilton-Jacobi equations. In the Section \ref{sec-DI} we recall some basics on the Optimal transport of probability measures and in particular about the equations of displacement interpolation. The final Section is devoted to prove the main results of the paper.\\

{\bf Acknowledgments}: We are grateful to Alberto Parmeggiani and Thierry Paul for the useful discussions and works on toroidal Pseudodifferential Opera\-tors and time propagation of semiclassical measures on the torus.

\section{Semiclassical measures}
\label{SEC-2}
Let us consider the flat torus $\Bbb T^n := (\Bbb R / 2\pi \Bbb Z)^n$.  
The  class of   symbols $b \in S^m_{\rho, \delta} (\mathbb{T}^n \times \mathbb{R}^n)$, $m \in \mathbb{R}$, $0 \le \delta$, $\rho \le1$, consist of those functions  
in $C^\infty (\mathbb{T}^n \times \mathbb{R}^n;\Bbb R)$ which are  $2\pi$-periodic in $x$ (that is, in each variable $x_j$, $1\leq j\leq n$) and for which
for all $\alpha, \beta \in \mathbb{Z}_+^n$ there exists $C_{\alpha \beta} >0$ such that $\forall$ $(x,\eta) \in \mathbb{T}^n \times \mathbb{R}^n$
\begin{equation}
\label{symb00}
|  \partial_x^\beta \partial_\eta^\alpha  b (x,\eta)  |   \le  C_{\alpha \beta m} \langle \eta \rangle^{m- \rho |\alpha| + \delta |\beta|}
\end{equation}
where $\langle\eta\rangle:=(1+|\eta|^2)^{1/2}$. In particular,  the set $S^m_{1,0} (\mathbb{T}^n \times \mathbb{R}^n)$ is denoted by $S^m (\mathbb{T}^n \times \mathbb{R}^n)$. 
The  toroidal Pseudodifferential Operator reads 
\begin{equation}
b(X,D) \psi(x):=(2\pi)^{-n}\sum_{\kappa \in\mathbb{Z}^n}\int_{\mathbb{T}^n}e^{i\langle x-y,\kappa \rangle}b(x,\kappa)\psi(y)dy, \quad \psi \in C^\infty (\mathbb{T}^n;\Bbb C), 
\end{equation}
see  \cite{R-T}. In particular, notice that it is given a map $b(X,D) : C^\infty (\mathbb{T}^n) \longrightarrow \mathcal{D}^\prime  (\mathbb{T}^n)$.  We recall that  $u \in \mathcal{D}^\prime  (\mathbb{T}^n)$ are the linear maps $u: C^\infty (\mathbb{T}^n) \longrightarrow \Bbb C$ such that $\exists$ $C>0$ and $k \in \Bbb N$, for which $|u(\phi)| \le C \sum_{|\alpha|\le k} \| \partial_x^\alpha \phi  \|_\infty$ $\forall \phi \in C^\infty (\Bbb T^n)$. 
Given a symbol $b\in S^m(\mathbb{T}^n\times\mathbb{R}^n)$, the (toroidal) Weyl quantization reads
\begin{equation}
\label{weyl}
\mathrm{Op}^w_\hbar(b)\psi(x) := (2\pi)^{-n}\sum_{\kappa \in\mathbb{Z}^n}\int_{\mathbb{T}^n}e^{i\langle x-y,\kappa \rangle}b(y,\hbar\kappa/2)\psi(2y-x)dy,\,\,\,\, \psi\in C^\infty(\mathbb{T}^n).
\end{equation}
In particular, it  holds  
\begin{equation}
\label{eq-O}
\mathrm{Op}^w_{\hbar} (b)  \psi (x) = (\sigma (X,D) \circ T_x \, \psi )(x)
\end{equation}
where  $T_x : C^\infty (\mathbb{T}^n) \rightarrow C^\infty (\mathbb{T}^n)$ defined as $(T_x \psi) (y) := \psi (2y-x)$  is linear, invertible and $L^2$-norm preserving, and $\sigma$ is a suitable toroidal simbol related to $b$, i.e.  $\sigma \sim \sum_{\alpha\geq 0}\frac{1}{\alpha!}\triangle_\eta^\alpha D_y^{(\alpha)} b(y,\hbar \eta / 2)\bigl|_{y=x}$, see Th. 4.2 in \cite{R-T} or also Th. 2.1 in \cite{P-Z}.\\ 
We say that a positive Radon measure with finite mass $\omega \in \mathcal{M}^{+} (\Bbb T^n \times \Bbb R^n)$ is a semiclassical measure associated with $\psi_\hbar \in L^2 (\mathbb{T}^n)$, $\| \psi_\hbar \|_{L^2} \le 1$  if there exists $\hbar_j \to 0^+$ as $j \to + \infty$ such that
\begin{equation}
\label{sem-def}
\lim_{j \to +\infty}
\langle \psi_{\hbar_j}, \mathrm{Op}^w_{\hbar_j} (\phi)  \psi_{\hbar_j} \rangle_{L^2} = \int_{\mathbb{T}^n \times \mathbb{R}^n} \phi (x,\xi) d\omega (x,\xi)
\end{equation}
for any test function $\phi \in C^\infty (\Bbb T^n \times \Bbb R^n)$ satisfying the phase space Fourier representation (see \cite{G-P})
\begin{equation}
\label{pf-T}
\phi (x,\xi) = (2\pi)^{-n} \int_{\mathbb{R}^n} \sum_{q \in \mathbb{Z}^n}    \widehat{\phi} (q,p) e^{ i (\langle p , \xi \rangle + \langle q , x \rangle)}   dp
\end{equation}
for some compactly supported  $\widehat{\phi} : \mathbb{Z}^n \times \mathbb{R}^n \to \Bbb R$, see Section 2.1.3 in \cite{TP-Z}. 

\section{A quick overview of weak KAM theory}
\label{wKAM}
The  weak KAM theory deals with a class of  Lipschitz continuous solutions of the Hamilton-Jacobi equation
\begin{equation}
\label{def-eff0}
H (x,\nabla_x v(x))  =  c[0]
\end{equation}
in the general assumption of Tonelli Hamiltonians  $H \in C^\infty (\mathbb{T}^n \times \mathbb{R}^n;\mathbb{R})$, that is to say, for functions $H$ such that $\eta \mapsto H(x,\eta)$ is strictly convex and uniformly superlinear in the fibers of the canonical projection $\pi : \mathbb{T}^n \times \mathbb{R}^n \longrightarrow \mathbb{T}^n$.  
The value $c[0]$ is called the {\it critical value} for which there exist solutions, and  it can be expressed by the inf-sup formula
\begin{equation}
\label{def-eff}
c[0]  \ =  \inf_{v \in C^\infty  (\mathbb{T}^n;\mathbb{R})} \  \sup_{x \in \mathbb{T}^n}  \ H(x,\nabla_x v(x)) 
\end{equation}
see for example \cite{E1}. If $H =  |p|^2 / 2m + V(x)$ then 
\begin{equation}
c[0] = \max_{y \in \Bbb T^n} V(y).
\end{equation}
The Lax-Oleinik semigroup of positive and negative type is defined as
$$
T_t^{\mp} u (x) :=  \inf_{\gamma} \Big\{ u(\gamma(0))  \pm  \int_0^t L(\gamma(s),\dot{\gamma}(s))  \  ds \Big\}, \quad u \in C^{0,1} (\Bbb T^n; \Bbb R), 
$$ 
where the infimum is taken over all continuous piecewise $C^1$ curves $\gamma : [0,1] \to \Bbb T^n$ such that $\gamma(t)=x$.  In particular, by defining $A^{0,t}(\gamma) := \int_0^t L(\gamma (\tau), \dot{\gamma}(\tau)) d  \tau$
\begin{equation} 
h_t (y,x) :=  \inf_{\gamma}  A^{0,t}(\gamma)
\end{equation} 
with $\gamma (0) = y$ and $\gamma(t) = x$, one can prove (see for example Prop. 4.1 in \cite{F-F}) that $h_t$ is continuous. Furthermore, it follows that 
$$
T_t^{-} u (x) =  \min_{y \in \Bbb T^n} \left\{ u(y) +  h_t (y,x)  \right\}, \quad T_t^{+} u (x) =  \max_{y \in \Bbb T^n} \left\{ u(y) -  h_t (x,y)  \right\}. 
$$
A  function $S_{-} \in C^{0,1}( \mathbb{T}^n ; \Bbb R)$   is said to be a {\it weak KAM solution of negative type} for (\ref{def-eff0}) if $\forall$ 
$t\ge0$
\begin{equation}
\label{back-}
T_t^{-} S_{-}  = S_{-}  - t \, c[0] ,  
\end{equation}
whereas  it is said to be a {\it weak KAM solution of positve type} if  $\forall$ $t\ge0$ 
\begin{equation}
\label{back+}
T_t^{+} S_{+}  = S_{+} + t \, c[0],
\end{equation}
see Def. 4.7.6 in \cite{F}. 
For any weak KAM solution it holds 
\begin{equation}
\label{inc-G}
\overline{{\rm Graph} (\nabla_x S_\pm )}  \subset  \{ (x,\eta) \in \Bbb T^n \times \Bbb R^n \ | \ H(x,\eta) = c [0] \}.
\end{equation}
Furthermore, the graphs are invariant under the backward (resp. forward) Hamiltonian flow, namely 
\begin{eqnarray}
\phi_H^t ( {\rm Graph}( \nabla_x S_{-}  ) ) &\subseteq& {\rm Graph}( \nabla_x S_{-} ) \quad \forall t \le 0
\\
\phi_H^t ( {\rm Graph}( \nabla_x S_+ ) ) &\subseteq& {\rm Graph}( \nabla_x S_+ ) \quad \forall t \ge 0
\label{forward-incl}
\end{eqnarray}
see Theorems 4.9.2 and 4.9.3 in \cite{F}. Moreover, it is proved that the maps $x \longmapsto (x, \nabla_x S_\pm)$ are continuous on  ${\rm  dom}(\nabla_x S_\pm) := \{ x \in \Bbb T^n \ | \  \exists \ \nabla_x S_{\pm} (x) \}$.
As showed within Th. 7.6.2 of \cite{F}, all the Lipschitz continuous weak KAM solutions of negative type coincide with the so-called {\it viscosity solutions} in the sense of \cite{C-L}.  

\section{The equations of displacement interpolation}
\label{sec-DI}
Let $X,Y$ be sets and $c : X \times Y \longrightarrow (-\infty,+\infty]$. A function $\psi : X \longrightarrow \Bbb R \cup \{ + \infty \}$ is said to be $c$-convex if it is not identically $+ \infty$ and there exists $\zeta : Y \longrightarrow \Bbb R \cup \{ \pm \infty \}$ such that 
$
\psi(x) = \sup_{y \in Y} (\zeta(y) - c(x,y))$ $\forall x \in X$.
Let $L \in C^2 (\Bbb T^n \times \Bbb R^n)$ be a Tonelli Lagrangian, and 
$A^{0,1}(\gamma) := \int_0^1 L(\gamma (\tau), \dot{\gamma}(\tau)) d  \tau$ the related  Lagrangian Action. Define the cost function 
\begin{equation} 
c^{0,1} (x,y):= \inf_\gamma A^{0,1}(\gamma)  
\end{equation} 
over all continuous piecewise $C^1$ curves $\gamma : [0,1] \to \Bbb T^n$ such that $\gamma(0)=x$ and  $\gamma (1) = y$. The related optimal transport cost reads 
\begin{equation} 
C^{0,1}(\mu,\nu) := \inf_{\pi \in \Pi (\mu,\nu)} \int_{\Bbb T^n \times \Bbb T^n}  c^{0,1} (x,y) \, d\pi(x,y).
\end{equation}   
Let $\sigma_0, \sigma_1 \in \mathcal{P}(\Bbb T^n)$ be such that $C^{0,1}(\sigma_0,\sigma_1) < + \infty$.  
Let   $\{\sigma_t \}_{0 \le t \le 1}  \in \mathcal{P}(\Bbb T^n)$ be a displacement interpolation of $\sigma_0$ and $\sigma_1$ with respect to the Lagrangian Action 
$A^{0,1}(\gamma)$. More precisely, the path $\{\sigma_t \}_{0 \le t \le 1}$ is linked to a minimizing curve for 
\begin{equation} 
\label{Min-A}
\inf_{\gamma} \int_{\Omega}    \ \int_0^1 L( \gamma (\tau,\omega), \partial_t \gamma (\tau,\omega)) d  \tau \ \rm{d} \mathbb{P}(\omega)
\end{equation}
where the infimum is over all the random curves $\gamma : [0,1] \times \Omega \longrightarrow \Bbb T^n$ such that ${\rm Law} (\gamma(\tau,\cdot \,) ) = \sigma_\tau$ for $0 \le \tau \le 1$; see Theorem 7.21 in \cite{V}. Then,  the following {\it equations of displacement interpolation} are fulfilled
\begin{itemize} 
\item[{\bf a.}]
$\partial_t \sigma_t (x) + {\rm div}_x (\xi (t,x)  \sigma_t (x) ) = 0$
\item[{\bf b.}]
$\nabla_v L(x,\xi (t,x)) = \nabla_x u (t,x)$ 
\item[{\bf c.}]
$\partial_t  u(t,x)   +   H(x,\nabla_x u (t,x)) =  0$   \qquad  $u(0,\cdot \, )$  is  $c$ - convex,  
\end{itemize} 
where the cost fuction is $c = c^{0,1} (x,y)$ as briefly outlined in chapter 13 of \cite{V}. In this setting, the vector field in the continuity equation can be equivalently written as $\xi (t,x) = \nabla_p H (x, \nabla_x u (t,x))$. 

\begin{remark}
\label{REM-convex}
In our paper $H = |p|^2 /2m + V(x)$ and hence $\xi (t,x) = \nabla_x u(t,x) / m$. Furthermore, our initial data is $u(0,\cdot \, ) =  S_+$, namely a  weak KAM solution of positive type for  the stationary Hamilton-Jacobi equation. Whence, $\forall t \ge 0$
\begin{equation}
T_t^{+} S_{+} -  t \, c[0]  = S_{+} 
 \end{equation}
namely 
\begin{equation}
\max_{y \in \Bbb T^n} \left\{ S_+ (y) -  h_t (x,y)  \right\} -  t \, c[0]  = S_{+} (x)
\end{equation}
which reads in the c - convex condition 
\begin{equation}
\max_{y \in \Bbb T^n} \left\{ S_+ (y) -  t \, c[0]  -  h_t (x,y) \right\}   = S_{+} (x).
\end{equation}
For $t = 1$ the cost function is in fact $h_1 (x,y) = c^{0,1} (x,y)$. Now, easily see that the function $S_+ (x) -  t \, c[0]$ is a solution of the equation ${\bf c}$ (equivalently $S_+$ solves the stationary H-J) and that the related continuity equation
$$
\partial_t \sigma_t (x) + {\rm div}_x \Big( \frac{1}{m} \nabla_x S_+ (t,x)  \sigma_t (x) \Big) = 0
$$
is solved by  $\sigma_t = (\Psi^t)_\sharp (\sigma_0)$ with $\Psi^t (x) := \pi \circ \phi_{H}^t (x,\nabla_x S_+ (x))$, as shown in Lemma \ref{CD}. 
We now recall the equivalence between optimal transport problems (i) - (iii) in Theorem 7.21 of \cite{V}, namely the link between optimal transference plans and displacement interpolations with respect to the Lagrangian Action. Moreover,    in view Theorem 12 - Proposition 1 of \cite{B-B} about Kantorovich optimal pairs (in our paper 
$(S_+, S_-)$)  we can apply Theorem 4.2 of \cite{F-F}  in the assumption $\sigma_0 \in \mathcal{P}_{ac} (\Bbb T^n)$. Thus, the path of measures $\sigma_t = (\Psi^t)_\sharp (\sigma_0)$ is a displacement interpolation in the sense of (\ref{cost00}). 
\end{remark}

\section{Main results}
\label{Results}

{\bf Proof of Theorem 1} Thanks to the setting of $\varphi_\hbar$,  any semiclassical measure $\omega_0 \in \mathcal{M}^+ (\Bbb T^n \times \Bbb R^n)$ associated with $\varphi_\hbar$ given by (\ref{def-wave}) takes the form 
\begin{equation} 
\omega_0 (x,p) = \delta ( p - \nabla_x S_+ (x) ) \sigma_0 (x)
\end{equation}
where $\sigma_0 \in \mathcal{P}_{ac} (\Bbb T^n)$, see Remark \ref{Rem-sem}. Hence, $\omega_0 \in \mathcal{P} (\Bbb T^n \times \Bbb R^n)$.\\
Let $H:= \frac{1}{2m} |p|^2 + V(x)$ and $\phi_{H}^t : \Bbb T^n \times \Bbb R^n \to \Bbb T^n \times \Bbb R^n$ the Hamiltonian flow. Applying Lemma \ref{AB}, the push forward   $\omega_t := (\phi_{H}^t)_\sharp \omega_0 \in \mathcal{P} (\Bbb T^n \times \Bbb R^n)$  reads for $t \ge 0$
\begin{equation} 
\omega_t (x,p) = \delta ( p - \nabla_x S_+ (x) ) \sigma_t (x)
\end{equation} 
where  $\sigma_t \in \mathcal{P} (\Bbb T^n)$. Thanks to Lemma \ref{CD},  this is a distributional solution for
\begin{equation}
\partial_t \sigma_t (x) +  {\rm div}_x  \Big( \frac{1}{m} \nabla_x S^+ (x)  \sigma_t (x) \Big) = 0
\end{equation}
which is fulfilled also by a continuous representative in the sense of Lemma 8.2.1 shown in \cite{A}.\\  
Any semiclassical limit ${\rm w}$ for the Wigner transform of $\psi (t,x) := (U_\hbar (t)  \varphi_\hbar) (x)$ in  $L^\infty ([0,1] ; A^\prime)$ solves the Liouville equation in the distributional sense 
\begin{equation}
\label{Liou-99}
\int_0^1 \int_{\mathbb{T}^n \times \Bbb R^n} [ \partial_s f (s,x,p) +   \{ H , f  \} (s,x,p) ]  d{\rm w}_s (x,p) ds  = 0 
\end{equation}
$\forall f \in C^\infty_c ((0,1) \times \Bbb T^n \times \Bbb R^n;\Bbb R)$ and moreover it holds the additional regularity $C([0,1];  \mathcal{P} (\Bbb T^n \times \Bbb R^n))$, as shown by  Theorem 4.1 - Remark 4.2 in \cite{TP-Z}. To conclude, the Liouville equation (\ref{Liou-99}) is linked to a smooth vector field $(p,-\nabla_x V(x))$
and hence it holds the uniqueness for the solutions in  $C([0,1];  \mathcal{P} (\Bbb T^n \times \Bbb R^n))$ which gives ${\rm w}_t = (\phi_{H}^t)_\sharp \omega_0$.   This implies,  the equality ${\rm w}_t = \omega_t$ for all $t \ge 0$. By recalling Remark \ref{REM-convex}  we conclude that  $\pi_\sharp (\omega_t) = \sigma_t$ for a.e. $0 \le t \le 1$ equals a path of continuous displacement interpolations in the sense of (\ref{cost00}) and (\ref{eqTR01}).    $\Box$  

\begin{remark}
\label{Rem-sem}
To prove that any semiclassical measure $\omega_0 \in \mathcal{M}^+ (\Bbb T^n \times \Bbb R^n)$, in the sense of (\ref{sem-def}), associated with $\varphi_\hbar$ as in (\ref{def-wave}) takes the form 
\begin{equation} 
\omega_0 (x,p) = \delta ( p - \nabla_x S_+ (x) ) \sigma_0 (x)
\end{equation}
we apply the same arguments shown in  Theorem  4.9 of \cite{TP-Z}. In fact, the proof is  based on the application of the following properties which recover the ones assumed in the present paper: 
\begin{itemize}
\item[i.] $a_\hbar \in H^1 (\Bbb T^n;\Bbb R)$ where $\| a_\hbar \|_{L^2} = 1$, $\| \hbar \nabla a_\hbar \|_{L^2} \to 0$ as $\hbar \to 0$, $a^2 (x) dx \rightharpoonup \sigma_0$ weakly as measures on $\Bbb T^n$,
\item[ii.]  ${\rm supp}(\sigma_0) \subseteq {\rm dom}(\nabla S_+)$,
\item[iii.]  $S_+ : \Bbb T^n \to \Bbb R$ is Lipschitz continuous,
\item[iv.]  $x \mapsto \nabla_x S_+ (x)$ is continuous on  ${\rm dom}(\nabla S_+)$.
\end{itemize}
The main difference here is that  we are dealing with a less general class of Hamilton-Jacobi equations (i.e. when $P=0$) and furthermore we are not assuming the absolute continuity  $\sigma_0 \ll \pi_\sharp (\mu_P)$ where $\mu_P$ is some invariant and Action minimizing measure.  In our paper we additionally assume that $\sigma_0 \ll  \mathcal{L}^n$ but this is not necessary for the proof of the semiclassical convergence  to the monokinetic measures $\omega_0$.  
\end{remark}

\begin{lemma}
\label{AB}
Let $S_+ : \Bbb T^n \to \Bbb R$ a Lipschitz continuous weak KAM solution of positive type for the H-J equation $\frac{1}{2m} |\nabla_x S_+ (x)|^2 + V(x) = \max_{y \in \Bbb T^n} V(y)$ and for some $\sigma_0 \in \mathcal{P} (\Bbb T^n)$ assume ${\rm supp}(\sigma_0) \subseteq{\rm  dom}(\nabla_x S_+)$.  Define $\omega_0 (x,p) := \delta ( p - \nabla_x S_+ (x) ) \sigma_0 (x)$. Let $H:= \frac{1}{2m} |p|^2 + V(x)$ and denote by $\phi_{H}^t : \Bbb T^n \times \Bbb R^n \to \Bbb T^n \times \Bbb R^n$ the Hamiltonian flow. Then, the push forward   $\omega_t := (\phi_{H}^t)_\sharp \omega_0 \in \mathcal{P} (\Bbb T^n \times \Bbb R^n)$  reads for $t \ge 0$
\begin{equation} 
\omega_t (x,p) = \delta ( p - \nabla_x S_+ (x) ) \sigma_t (x)
\end{equation} 
with  $\sigma_t \in \mathcal{P} (\Bbb T^n)$ as in (\ref{Lab}).
\end{lemma}
\proof
For any test function $f \in C^\infty_c (\Bbb T^n \times \Bbb R^n;\Bbb R)$ it holds
\begin{equation} 
\label{omega56}
\int_{\Bbb T^n \times \Bbb R^n}  f (x,p) \ d \omega_t (x,p)  =  \int_{\Bbb T^n \times \Bbb R^n}  f \circ \phi_H^{t} (x,p) \ d \omega_0 (x,p) 
\end{equation} 
and, by the assumption on $\omega_0$, 
\begin{equation} 
\label{34-Int}
\int_{\Bbb T^n \times \Bbb R^n}  f (x,p) \ d \omega_t (x,p)  =  \int_{\Bbb T^n}  f \circ \phi_H^{t} (x,\nabla_x S_+ (x) ) \ d \sigma_0 (x). 
\end{equation} 
Indeed,  we recall that the map $x \longmapsto (x, \nabla_x S_+)$ is continuous when restricted on the set ${\rm  dom}(\nabla_x S_+) := \{ x \in \Bbb T^n \ | \  \exists \ \nabla_x S_+ (x) \}$ which is a Borel set. 
 Thus,  
\begin{equation} 
x \mapsto  f \circ \phi_H^{t} (x,\nabla_x S_+ (x) )
\end{equation} 
is  a continuous map on ${\rm  dom}(\nabla_x S_+)$ and hence also on ${\rm supp}(\sigma_0)$. Whence, the integral (\ref{34-Int}) is well posed.
Furthermore, remind that
\begin{equation} 
\phi_H^{t} ( {\rm Graph}( \nabla_x S_+ ) ) \subseteq {\rm Graph}( \nabla_x S_+ )
\end{equation} 
for any $t \ge 0$. 
Thus, for $\Psi^{t} (x) := \pi \circ \phi_H^{t} (x,\nabla_x S_+ (x) )$  and $t \ge 0$
\begin{equation} 
\label{eq-37phi}
\phi_H^{t} (x,\nabla_x S_+ (x) ) =  ( \Psi^{t} (x), \nabla_x S_+ ( \Psi^{t} (x)) ) 
\end{equation} 
In particular, any map $\Psi^{t} : {\rm  dom}(\nabla_x S_+ ) \to \Psi^{t} ({\rm  dom}(\nabla_x S_+)) \subseteq {\rm  dom}(\nabla_x S_+)$ is continuous and one to one. In addition, notice that $t \mapsto \Psi^{t} (x)$ is absolutely continuous for any $x \in {\rm  dom}(\nabla_x S_+)$. 
To conclude, thanks to (\ref{eq-37phi})   the integral (\ref{34-Int}) can be rewritten as 
\begin{equation} 
\int_{\Bbb T^n}  f  ( \Psi^{t} (x)   , \nabla_x S_+ ( \Psi^{t} (x)   ) )   \ d \sigma_0 (x). 
\end{equation} 
By defining 
\begin{equation} 
\label{Lab}
\sigma_t  :=  (\Psi^{t} )_\sharp  \sigma_0
\end{equation} 
we recover for (\ref{omega56}) the form
\begin{eqnarray} 
\int_{\Bbb T^n}  f  ( x , \nabla_x S_+ (x) ) \   d\sigma_t (x).
\end{eqnarray} 
$\Box$ 


\endproof

\begin{lemma}
\label{CD}
Let $\sigma_t  :=   (\Psi^{t})_\sharp  (\sigma_0)$ be as in (\ref{Lab}). Then, 
\begin{equation}
\int_0^1 \int_{\Bbb T^n}   |\nabla_x S^+ (x)|  \, d \sigma_t (x)  < + \infty
\end{equation}
and $\forall f \in C^\infty_c ( (0,1) \times \Bbb T^n;\Bbb R)$
\begin{equation}
\label{eq-L}
\int_0^1 \int_{\Bbb T^n} \Big( \partial_t f (t,x) +  \nabla_x f(t,x) \cdot  \frac{1}{m} \nabla_x S^+ (x) \Big) d\sigma_t (x)  dt   = 0.
\end{equation}
Morever, there exists a narrowly continuous curve $t \in [0,1] \to \widetilde{\sigma}_t \in \mathcal{P} (\Bbb T^n)$ such that $\sigma_t =  \widetilde{\sigma}_t$ for $\mathcal{L}^1$ - a.e. $t \in(0,1)$.   
\end{lemma}
\proof 
About the first condition, we recall the setting of $\sigma_t$ and the assumption ${\rm supp} (\sigma_0) \subseteq {\rm dom}(\nabla_x S_+)$,   
\begin{equation}
\int_0^1 \int_{\Bbb T^n}   |\nabla_x S^+ (x)|  \, d \sigma_t (x)   =   \int_0^1 \int_{\Bbb T^n}   |\nabla_x S^+ (\Psi^{t} (x))|  \, d \sigma_0 (x) . 
\end{equation}
In particular, recalling (\ref{inc-G}) and (\ref{forward-incl}),  it follows directly  
\begin{equation}
\sup_{x \in {\rm supp} (\sigma_0)}  |\nabla_x S^+ (\Psi^{t} (x))| \le \sup_{y \in  {\rm dom}(\nabla_x S_+) }  |\nabla_x S^+ (y)| < + \infty.
\end{equation}
Furthermore,  the integral in (\ref{eq-L}) reads 
\begin{equation}
\int_0^1 \int_{\Bbb T^n} \Big( \partial_t f (t,\Psi^{t} (x)) +  \nabla_x f(t,\Psi^{t} (x)) \cdot  \frac{1}{m} \nabla_x S^+ (\Psi^{t} (x)) \Big) d\sigma_0 (x)  dt   
\end{equation}
and recalling the setting of $\Psi_t (x)$
\begin{equation}
\int_0^1 \int_{\Bbb T^n} \Big( \partial_t f (t,\Psi^{t} (x)) +  \nabla_x f(t,\Psi^{t} (x)) \cdot  \frac{d}{dt} \Psi^{t} (x)  \Big) d\sigma_0 (x)  dt .
\end{equation}
This expression reads equivalently
\begin{eqnarray}
\int_0^1 \int_{\Bbb T^n}  \frac{d}{dt}   f (t,\Psi^{t} (x)) \, d\sigma_0 (x)  dt   &=&  \int_{\Bbb T^n} \int_0^1 \frac{d}{dt}   f (t,\Psi^{t} (x)) dt   d\sigma_0 (x)
\\
&=&  \int_{\Bbb T^n}   \tilde{f} (t,\Psi^{t} (x)) |_{0}^1  \, d\sigma_0 (x).
\end{eqnarray}
To conclude,  we notice that any test functions $f \in C^\infty_c ( (0,1) \times \Bbb T^n;\Bbb R)$ the (vanishing) smooth extention $\tilde{f}$ at $t=0$ and $t=1$ fulfills $ \tilde{f} (t=0,\Psi^{t=0} (x)) = \tilde{f} (t=1,\Psi^{t=1} (x)) = 0$. 
By applying Lemma 8.1.2 of \cite{A} it follows the existence of a narrowly continuous curve $t \in [0,1] \to \widetilde{\sigma}_t \in \mathcal{P} (\Bbb T^n)$ such that $\sigma_t =  \widetilde{\sigma}_t$ for $\mathcal{L}^1$ - a.e. $t \in(0,1)$.
$\Box$

\begin{remark}
Working with narrowly continuous curves $t \in [0,1] \to \widetilde{\sigma}_t \in \mathcal{P} (\Bbb T^n)$ means that 
$\widetilde{\sigma} \in C([0,1];  \mathcal{P} (\Bbb T^n) )$ and $\mathcal{P} (\Bbb T^n)$ is equipped with the L\'evy-Prokhorov distance of measures which metrizes the weak convergence. 
\end{remark}

\bigskip
\noindent
{\bf Proof of Theorem 2} Fix $\sigma_0 \in \mathcal{P}_{ac} (\Bbb T^n)$.  In view of Proposition 4.6 in \cite{P-Z} and Theorem 4.9 in \cite{TP-Z} there exists $\varphi_\hbar$ as in (\ref{def-wave}) such that $\omega_0 = \delta ( p - \nabla_x S_+) \sigma_0$ is the unique related semiclassical measure.\\ 
We suppose that the continuity equation 
\begin{equation}
\int_0^1 \int_{\Bbb T^n} 
\Big( \partial_t f(t,x) +   \nabla_x f(t,x) \cdot \frac{1}{m} \nabla_x S_+  (x)  d\sigma_t (x) dt \Big)
= 0
\end{equation}
with $f \in C^\infty_c ( (0,1) \times \Bbb T^n;\Bbb R)$ has a unique solution in  $C([0,1];  \mathcal{P} (\Bbb T^n) )$. Hence, this solution must coincide with a continuous representative $\widetilde{\sigma}_t$ of $\sigma_t  :=   (\Psi^{t})_\sharp  (\sigma_0)$ as in (\ref{Lab}).\\ 
Define the cotangent bundle lift 
$\widehat{\omega}_t := \delta ( p - \nabla_x S_+)\widetilde{\sigma}_t$. 
In particular, since the map $x \mapsto \nabla_x S_+ (x)$ is continuous on its domain, the lift $\widehat{\omega}_t$ fulfills $\widehat{\omega} \in C([0,1];  \mathcal{P} (\Bbb T^n \times \Bbb R^n))$. Recalling Lemma \ref{AB}, any $\widehat{\omega}_t$ equals   $\omega_t := (\phi_{H}^t)_\sharp \omega_0$. 
Solutions $\omega_t$ of the Liouville equation with the class of test functions $f \in C^\infty_c ( (0,1) \times \Bbb T^n \times \Bbb R^n;\Bbb R)$ are unique in $C([0,1];  \mathcal{P} (\Bbb T^n \times \Bbb R^n))$. 
Recalling  Theorem 4.1 and Remark 4.2 in \cite{TP-Z}, it follows that the solution $\omega_t$ coincides with the continuous path of semiclassical measures linked to the solution of the Schr\"odinger equation with our class of initial data $\varphi_\hbar$. 
$\Box$ 

\endproof

%




\end{document}